\documentclass{amsart}
\usepackage{quiver}
\usepackage[utf8]{inputenc}
\usepackage{amsmath,amssymb,amsthm,units, stmaryrd,stackrel,relsize,bm}
\usepackage{enumitem}
\usepackage{mathdots}
\usepackage{bussproofs}
\usepackage{mathtools}
\usepackage{adjustbox}
\usepackage{tikz,tikz-cd}
\usetikzlibrary{calc,decorations.pathmorphing}
\usetikzlibrary{arrows.meta}
\usetikzlibrary{fit,patterns,decorations.markings,matrix,3d}

\newsavebox{\pullback}
\sbox\pullback{%
\begin{tikzpicture}%
\draw (0,0) -- (1ex,0ex);%
\draw (1ex,0ex) -- (1ex,1ex);%
\end{tikzpicture}}
\newtheorem{theorem}{Theorem}
\newtheorem{theoremletter}{Theorem}

\newtheorem{lemma}[theorem]{Lemma}

\newtheorem{claim}[theorem]{Claim}
\newtheorem{corollary}[theorem]{Corollary}

\theoremstyle{definition}
\newtheorem{definition}[theorem]{Definition}

\theoremstyle{remark}
\newtheorem{remark}[theorem]{Remark}

\usepackage[
    colorlinks=true, 
    citecolor=red,
    linkcolor=blue,
    urlcolor=blue]{hyperref}
\newcommand\zfc{{\mathsf{ZFC}}}

\newcommand\kp{{\mathsf{KP}}}

\newcommand\glp{\mathsf{GLP}}
\newcommand\ig{\mathfrak{I}}
\newcommand\ja{\mathsf{J}}
\newcommand\lin{\mathsf{lin}}

\makeatletter
\def\Ddots{\mathinner{\mkern1mu\raise\p@
\vbox{\kern7\p@\hbox{.}}\mkern2mu
\raise4\p@\hbox{.}\mkern2mu\raise7\p@\hbox{.}\mkern1mu}}
\makeatother
\begin{document}
\title{The Logic of Correct Models}
%\subjclass[2020]{03B45, 03E99}
\author{J. P. Aguilera}
\address{1. Institute of Discrete Mathematics and Geometry, Vienna University of Technology. Wiedner Hauptstra{\ss}e 8--10, 1040 Vienna, Austria.}
\address{2. Department of Mathematics, University of Ghent. Krijgslaan 281-S8, B9000 Ghent, Belgium}
\email{aguilera@logic.at}
\author{F. Pakhomov}
\address{1. Department of Mathematics, University of Ghent. Krijgslaan 281-S8, B9000 Ghent, Belgium}
\address{2. Steklov Mathematical Institute of the Russian Academy of Sciences. Ulitsa Gubkina 8, Moscow 117966, Russia.}
\email{fedor.pakhomov@ugent.be}
%\begin{document}

\begin{abstract}
For each $n\in\mathbb{N}$, let $[n]\phi$ mean ``the sentence $\phi$ is true in all  $\Sigma_{n+1}$-correct transitive sets.'' Assuming G\"odel's axiom $V = L$, we prove the following graded variant of Solovay's completeness theorem:  the set of formulas valid under this interpretation is precisely the set of theorems of the linear provability logic $\glp.3$. We also show that this result is not provable in $\zfc$, so the hypothesis $V = L$ cannot be removed. 

As part of the proof, we derive (in $\zfc$) the following purely modal-logical results which are of independent interest: the logic $\glp.3$ coincides with the logic of closed substitutions of $\glp$, and is the maximal non-degenerate, normal extension of $\glp$.
\end{abstract}
\date{February, 2024}
\subjclass[2020]{03B45, 03E99}
\keywords{Elementary submodel, provability logic, completeness theorem, modal logic}
\clearpage
\maketitle
\numberwithin{equation}{section}
\setcounter{tocdepth}{1}
%\tableofcontents

\section{Introduction}
The study of arithmetical and set-theoretic interpretations of modal logic goes back to Friedman's 35th problem \cite{Fr75}. Friedman asked for an axiomatization of the set of valid closed formulas built from expressions of the form ``$\phi$ is provable.'' 
The problem was solved by Boolos \cite{Boo75} and independently by van Benthem and by Bernardi and Montagna, whose solution was later extended by the completeness theorem of Solovay \cite{So76}. Solovay also axiomatized the set of valid formulas built from expressions of the form ``$\phi$ holds in all models of $\zfc$ of the form $V_\kappa$.'' 
In this article, we present a graded version of Solovay's theorem axiomatizing the set of valid formulas built from expressions of the form ``$\phi$ holds in all $\Sigma_{n+1}$-correct models of set theory,'' assuming that G\"odel's axiom $V = L$ holds. 
% This extends Solovay's theorem, due to the observation that all transitive models of $\zfc$ of the form $V_\kappa$ are $\Sigma_1$-correct by the L\"owenheim-Skolem theorem. Let us state the result more precisely.

We consider expressions built from propositional variables $p,q,r,\hdots$ by means of Boolean connectives $\wedge,\vee\to\lnot$ and infinitely many operators $[0], [1], [2], \hdots$. Denote this language by $\mathcal{L}_{[<\omega]}$.
If $\phi$ is any formula built this way, $\phi$ can be interpreted as a set-theoretic sentence as follows: a \textit{correct-model interpretation} 
is an assignment $\phi \mapsto \llbracket \phi \rrbracket$
 of set-theoretic sentences to $\mathcal{L}_{[<\omega]}$-formulas which commutes with Boolean connectives and such that for all $\phi$ and all $[n]$, the sentence $\llbracket [n]\phi \rrbracket$ is an appropriate formalization of ``$\llbracket \phi \rrbracket$ holds in all $\Sigma_{n+1}$-correct transitive sets.'' Here, a set $M$ is $\Sigma_{n+1}$-correct if 
$(M, \in) \prec_{\Sigma_{n+1}} (V,\in),$
 i.e., $M$ is a $\Sigma_{n+1}$-elementary substructure of the universe. 
 %All the results in this article also go through in a slightly stronger form obtained by interpreting $[n]\phi$ as ``$\llbracket\phi\rrbracket$ holds in all $\Sigma_n$-correct transitive models of $V = L$.''
 A formula $\phi$ is \textit{valid} in a set $M$ if $\llbracket \phi \rrbracket$ holds in $M$ for every correct-model interpretation $\llbracket \cdot \rrbracket$. A formula $\phi$ is
 \textit{valid under the correct-model interpretation} if $\phi$ is valid in every $\Sigma_1$-correct $M$. If one defines validity of $\phi$ as ``for all correct-model interpretations we have $\zfc \vdash \llbracket \phi\rrbracket$,'' then all results in this article remain valid, and provably so in the theory $\zfc + \text{con}(\zfc)$.
%Observe that if we assume $V = L$ holds, then all $\Sigma^1_2$-correct models satisfy $V = L$, so the hypothesis is only relevant for the interpretation of formulas of the form $[0]\phi$.

We shall prove that the set of formulas valid under this interpretation coincides with the  natural extension of Japaridze's \cite{Ja88} provability logic $\glp$ by the well-known linearity axiom $.3$, assuming $V = L$.\\

\begin{theoremletter}\label{TheoremMain}
Suppose that $\zfc$ and $V = L$ hold. Let $\phi$ be an $\mathcal{L}_{[<\omega]}$-formula. 
Then,
the following are equivalent:
\begin{enumerate}
\item\label{TheoremMain1} $\phi$ is valid under the correct-model interpretation, and
\item\label{TheoremMain2} $\phi$ is a theorem of $\glp.3$.
\end{enumerate} 
\end{theoremletter}
Strictly speaking, Theorem \ref{TheoremMain} should be regarded as a theorem schema, with a distinct theorem for each $\phi$. This is solely due to the fact that the language of set theory does not allow uniformly discussing $\Sigma_n$-correct models for all $n \in\mathbb{N}$. See \S\ref{SectProof} for a discussion on this and on some variants of the correct-model interpretation for which the theorem also holds as a single sentence.

The logic $\glp.3$ and others are defined in \S \ref{SectGLP} below. As part of the proof of Theorem \ref{TheoremMain}, we obtain the following purely modal-logical characterization which is of independent interest and does not require the assumption that $V = L$:
\begin{theoremletter}\label{TheoremMainGLP}
Let $\phi$ be an $\mathcal{L}_{[<\omega]}$-formula. Then, the following are equivalent:
\begin{enumerate}
\item[(2)] \label{TheoremMainGLP1} $\phi$ is a theorem of $\glp.3$, 
\item[(3)] \label{TheoremMainGLP2} every closed substitution of $\phi$ is a theorem of $\glp$,
\item[(4)] \label{TheoremMainGLP3} $\phi$ is a theorem of some non-degenerate, normal extension of $\glp$.
\end{enumerate}
\end{theoremletter}
Theorem \ref{TheoremMainGLP} asserts that $\glp.3$ is the maximal non-degenerate, normal extension of $\glp$. The equivalence between (2) and (3) strengthens a result of Icard and Joosten \cite{IcJo12} who proved the equivalence for formulas in which $[n]$ occurs for at most one $n$. 

The proof of Theorems \ref{TheoremMain} and \ref{TheoremMainGLP} is included in \S \ref{SectProof}.
It follows from the proof of Theorem \ref{TheoremMainGLP} that the implication from \eqref{TheoremMain1} to \eqref{TheoremMain2} in the statement of Theorem \ref{TheoremMain} does not require the assumption that $V = L$; however, the converse direction does:
\begin{theoremletter}\label{TheoremVL}
The statement that the correct-model interpretation of $\mathcal{L}_{[<\omega]}$ satisfies axiom $.3$ fails in a generic extension of $L$ and is thus independent of $\zfc$.
\end{theoremletter}
Theorem \ref{TheoremVL} is proved in a significantly stronger form in \S \ref{SectTheoremC}.\\

$\glp$ has been the subject of extensive study due to its interesting model theory and its applications to the study of formal provability and proof theory of arithmetic. 
While modal logics are commonly studied with the method of relational semantics, $\glp$ has no non-trivial Kripke frames and thus requires move involved semantics given in terms of polytopological spaces.
The construction of such spaces usually depends on set-theoretic hypotheses beyond Zermelo-Fraenkel set theory, such as Jensen's \cite{Je72} $\Box_\kappa$ principle (see Blass \cite{Bl90}), the Axiom of Choice (see Beklemishev and Gabelaia \cite{BG13}, Fern\'andez-Duque \cite{FD14}, Shamkanov \cite{Sh20}, and also \cite{AgFD17,Ag22}, which build on Abashidze \cite{Ab85} and Blass \cite{Bl90}), indescribable cardinals and hyperstationary reflection (see Bagaria \cite{Ba19}, Bagaria, Magidor, and Sakai \cite{BMS15}, and Brickhill and Welch \cite{BrWe23}) and the structure of the Mitchell order on normal measures (see Golshani and Zoghifard \cite{GZ24}). For an overview, see Beklemishev and Gabelaia \cite{BG14}.

The logic $\glp$ is complete for various arithmetical interpretations in first- and second-order arithmetic (see Japaridze \cite{Ja88}, Ignatiev \cite{Ig93}, Cord\'on-Franco, Fern\'andez-Duque, Joosten, and Lara-Mart\'in \cite{CFFDJoLM17}, and Fern\'andez-Duque and Joosten \cite{FDJo18}) which have applications to Ordinal Analysis and the study of arithmetical theories (see e.g., Beklemishev \cite{Be04} and also \cite{BP22}). The correct-model interpretation for the one-modality fragment of $\glp$ is a central tool in the proof of the non-linearity theorem of \cite{AgPac}.
 
The main technical ingredient in the proof of Theorem \ref{TheoremMainGLP} is a modal-logical result, Theorem \ref{TheoremCompletenessJlin} in \S\ref{SectHL}, which is of independent interest. It asserts that an extension of Beklemishev's \cite{Be10} logic $\ja$ is sound and complete for the class of so-called hereditarily linear $J_n$-trees, or \textit{$J_n$-lines}. These had been previously studied by Icard \cite{Ic08,Ic11}.
Our proof is a modified ``canonical-model'' construction, as is typical in modal logic. However, for various reasons, the construction is more involved  than usual and requires considering the interaction between infinite and finite versions of the canonical models.

\section{The logic $\glp$}\label{SectGLP}
\begin{definition}
Let $\mathcal{L}_{[<\omega]}$ be the language of propositional logic augmented with operators $[n]$ for $n\in\mathbb{N}$ and their duals $\langle n\rangle$. The logic $\glp$ is the logic in the language $\mathcal{L}_{[<\omega]}$ given by all propositional tautologies and the following axioms for each $n\in\mathbb{N}$
\begin{enumerate}
\item $[n](\phi\to\psi) \to ([n]\phi \to [n]\psi)$;
\item $[n]([n]\phi\to\phi) \to [n]\phi$;
\item $[n]\phi \to [n+1]\phi$;
\item $\langle n\rangle \phi \to [n+1]\langle n\rangle \phi$;
\end{enumerate}
and closed under the rules modus ponens (from $\phi \to \psi$ and $\phi$ derive $\psi$) and necessitation (from $\phi$ derive $[n]\phi$).
\end{definition}

In particular, $\glp$ is a normal modal logic. Here, recall that a modal logic is called \textit{normal} if it contains the axiom 
$\Box(\phi\to\psi)\to(\Box\phi\to\Box\psi)$
for each modality
and is closed under modus ponens and necessitation for each modality. 

Let $\mathcal{L}_{[<n]}$ be the language of propositional logic augmented with operators $[m]$ for $m< n$ and their duals $\langle m\rangle$.  
The logic $\glp_n$ is defined like $\glp$ but restricting to the language $\mathcal{L}_{[<n]}$. Below and throughout let us assume that the language of $\mathcal{L}_{[<\omega]}$ is additionally enhanced with a constant symbol $\top$. The interpretation of $\top$ is $1 = 1$. We abbreviate $\lnot\top$ by $\bot$.

A \textit{Kripke frame} for $\mathcal{L}_{[<n]}$ is a structure 
$(W, <_0, <_1, \hdots, <_n),$ where each $<_k$ is a binary relation on $W$, called the \textit{accessibility relation}. Kripke frames can be endowed with interpretation functions, turning them into $\mathcal{L}_{[<n]}$-models. These  are functions $\phi\mapsto \llbracket \phi\rrbracket$ mapping $\mathcal{L}_{[<n]}$-formulas to subsets of $W$ such that $\llbracket \bot \rrbracket = \varnothing$, which commute with Booleans, and such that
\[\forall k\leq n\,\Big(w \in \llbracket \langle k \rangle \phi \rrbracket \leftrightarrow \exists v\, w <_k v\, \big( v \in \llbracket \phi \rrbracket \big) \Big).\]
We write 
$w \Vdash \phi$
to mean that $\phi$ holds true at $w \in W$, and write 
$W, w\Vdash \phi$
to emphasize that this evaluation takes place in the model $W$. Note: we might abuse notation by identifying frames with models.

A formula in $\mathcal{L}_{[<\omega]}$ is \textit{closed} if it contains no propositional variables. Thus, it is built up from the symbol $\top$ using Boolean connectives and operators $[n]$. A \textit{closed substitution} of a formula $\phi$ is obtained by substituting a closed $\mathcal{L}_{[<\omega]}$-formula  for each propositional variable.

\begin{definition}
The logic $\glp.3$ is obtained by extending $\glp$ with the linearity axiom $.3$, for each $n\in\mathbb{N}$:
\[[n]([n]A\to B) \vee [n]([n]B\wedge B \to A).\]
\end{definition}
$\glp_n.3$ is defined analogously.

\begin{theorem}[$V = L$] \label{TheoremValidity}
Suppose $\phi$ is a theorem of $\glp.3$. Then $\phi$ is valid under the correct-model interpretation.
\end{theorem}
\proof
Here and throughout the article, we abuse notation by identifying formulas $\phi$ with their interpretations $\llbracket \phi \rrbracket$.

Clearly validity is preserved under modus ponens. For necessitation, suppose $\phi$ is valid in every $\Sigma_1$-correct set. Then, if $M$ is a $\Sigma_1$-correct set and $M \models ``N$ is $\Sigma_{n+1}$-correct,'' then $N$ is $\Sigma_1$-correct, so $\phi$ is valid in $N$. Thus, $[n]\phi$ is valid. Here, observe that $\Sigma_1$-correct sets are models of $\kp$ (with infinity), and this theory is enough to uniformly speak of partial satisfaction classes for $\Sigma_n$ formulas, which is needed to formalize $\Sigma_n$-correctness.
It remains to verify that the axioms are valid. The axioms 
$[n](\phi\to\psi) \to ([n]\phi \to [n]\psi)$
and 
$[n]\phi \to [n+1]\phi$
are clearly valid, as is the principle $[n]\phi \to [n][n]\phi$. Thus, the correct-model interpretation includes the logic $\mathsf{K4}$, which is enough to varry out the usual proof of L\"ob's theorem, yielding the validity of  
$[n]([n]\phi\to\phi) \to [n]\phi$. We remark that, although L\"ob's theorem is usually applied to provability, the same proof yields validity for the correct-model interpretation. %See \cite{AgPac} for a proof.
The validity of 
$\langle n\rangle \phi \to [n+1]\langle n\rangle \phi$
follows from the observation that the predicate ``$M$ is $\Sigma_{n+1}$-correct'' is $\Pi_{n+1}$. We have shown that $\glp$ is valid, without making use of the hypothesis $V = L$.

Before moving on, we make an observation. Suppose $M$ is $\Sigma_1$-correct. Suppose $x \in M$. By $V = L$ and $\Sigma_1$-correctness, we have 
$M \models \exists \alpha\, x \in L_\alpha,$
so $M \models V = L$. By G\"odel's condensation lemma, we have $M = L_\gamma$ for some $\gamma$.

Now, work in a model of the form $L_\xi$ and let $\phi$ and $\psi$ be two formulas. We show that 
\[[n] ([n]\phi\to\psi)\vee [n]([n]\psi\wedge\psi\to\phi)\]
holds.
Let $\alpha$ and $\beta$ be, respectively, least such that $L_\alpha$ and $L_\beta$ are $\Sigma_{n+1}$-correct, $L_\alpha\not\models\phi$, $L_\beta\not\models \psi$. If no such $\alpha$ exists, put $\alpha = \infty$, and similarly for $\beta$.\smallskip

\noindent\textsc{Case I}: $\alpha<\beta$. Let $M$ be $\Sigma_{n+1}$-correct. Then, $M$ is of the form $L_\gamma$. Suppose 
$M \models [n]\phi$.
Then $\gamma \leq\alpha$ for otherwise $L_\alpha \in L_\gamma$ and by correctness
\[M\models \text{``$L_\alpha$ is $\Sigma_{n+1}$-correct,''}\]
so that $M\models \lnot [n]\phi$. Then, $\gamma\leq\alpha <\beta$ and hence $L_\gamma\models\psi$ by choice of $\beta$. We have showed that 
$[n] ([n]\phi\to\psi)$
holds.\smallskip

\noindent\textsc{Case II}: $\beta\leq\alpha$. Let $M$ be $\Sigma_{n+1}$-correct. Then, $M$ is of the form $L_\gamma$. Suppose 
$M \models [n]\psi \wedge\psi.$
Then $\gamma < \beta$ for if $\gamma = \beta$ then $L_\gamma\not\models \psi$, and if $L_\beta \in L_\gamma$ then above we have $M\models \lnot [n]\psi$. Thus $\gamma <\beta \leq\alpha$ and hence $L_\gamma\models\phi$ by choice of $\alpha$. We have shown that 
$[n]([n]\psi\wedge\psi\to\phi)$
holds. This proves the theorem.
\endproof

\section{Hereditarily linear frames}\label{SectHL}
Towards the proofs of the main theorems, we take a technical detour in which we axiomatize the class of so-called hereditarily linear $J$-frames via an extension of Beklemishev's logic $\ja$.

\begin{definition}
The logic $\ja_n$ is the logic in the language $\mathcal{L}_{[<n]}$ given by all propositional tautologies and the following axioms:
\begin{enumerate}
\item $[k](\phi\to\psi) \to ([k]\phi \to [k]\psi)$;
\item $[k]([k]\phi\to\phi) \to [k]\phi$, for $k < n$;
\item $[k]\phi \to [k][m]\phi$, for $k \leq m < n$;
\item $[k]\phi \to [m][k]\phi$, for $k \leq m < n$;
\item $\langle k\rangle \phi \to [m]\langle k\rangle \phi$, for $k < m < n$;
\end{enumerate}
and closed under the rules modus ponens (from $\phi \to \psi$ and $\phi$ derive $\psi$) and necessitation (from $\phi$ derive $[k]\phi$). 
\end{definition}

\begin{definition}
The logic $\ja_n.\lin$ is obtained by extending $\ja_n$ by the following ``pseudo-linearity'' axioms, for each $m < n$:
\begin{align*}
[m] \bigg(\Big(\bigwedge_{m\leq k< n}[k] \phi \Big) \to \psi \bigg) 
\vee 
[m] \bigg(\Big(\bigwedge_{m\leq k< n}[k]\psi \wedge \psi \Big) \to \phi \bigg) 
\end{align*}
\end{definition}

\begin{lemma}
$\ja_n.\lin$ is a sublogic of $\glp_n.3$.
\end{lemma}
\proof
Immediate.
\endproof

\begin{definition}
A $J_n$-frame is a frame $(T, <_0, \hdots, <_{n-1})$ where each $<_i$ is irreflexive, transitive, and converse wellfounded, and such that the following three properties hold for $k<m$:
\begin{enumerate}
\item if $x <_m y <_k z$, then $x<_k z$;  
\item if $x <_k y <_m z$, then $x<_k z$; and
\item if $x<_k z$ and $y<_m z$, then $x<_k y$.
\end{enumerate}
\end{definition}

\begin{definition}
A $J_n$-frame $(T, <_0, \hdots, <_{n-1})$ is \textit{hereditarily linear} if for all $x,y\in T$ such that $x\neq y$, there is $k<n$ such that $x <_k y$ or $y <_k x$.
\end{definition}

\begin{lemma} \label{LemmaHLOrder}
Suppose $(T, <_0, \hdots, <_{n-1})$ is a hereditarily linear $J_n$-frame. 
Then, for all $x,y\in T$ such that $x\neq y$, there is precisely one $k$ such that one of $x <_k y$ or $y <_k x$ holds. For such $k$, precisely one of $x <_k y$ or $y <_k x$ holds.
\end{lemma}
\proof
Hereditary linearity yields some $k$ such that, say, $x <_k y$. Since $J_n$-frames are irreflextive, we have $y \not <_k x$.
The fact that $T$ is a $J_n$-frame implies that there can only be one such $k$, for  $x <_l y$ or $y<_l x$ would imply $x <_{\min\{k,l\}} x$; if $x <_l y$ this follows from propery (3) of $J_n$-frames, if $y<_lx$ and $l<k$ this follows from property (1), and if $y<_lx$ and $k<l$ this follows from property (2).
\endproof

\subsection{Properties of hereditarily linear frames}
In this section, we record some properties of hereditarily linear frames which will not be used later, but which might be useful to note.
\begin{definition}
A $J_n$-frame $(T, <_0, \hdots, <_{n-1})$ is \textit{stratified} if $x <_m y$ and $z <_k y$ imply $x<_m z$ whenever $k < m$.
\end{definition}

\begin{lemma} \label{LemmaStratified}
Suppose $(T, <_0, \hdots, <_{n-1})$ is a hereditarily linear $J_n$-frame. Then, $(T, <_0, \hdots, <_n)$ is stratified.
\end{lemma}
\proof
Suppose $x <_k z$ and $y <_m z$, where $k < m$. Since $T$ is a $J_n$-frame, this implies $x\neq y$.
By hereditary linearity, we have some $l$ such that $x <_l y$ or $y <_l x$. Suppose that $y <_l x$. If $l \leq k$, then we have $y <_l z$, since $T$ is a $J_n$-frame, which together with $y <_m z$ implies $z <_l z$ -- a contradiction. If $k < l$, then we have $y <_k z$ by Lemma \ref{LemmaJFrame}, which again by Lemma \ref{LemmaJFrame} contradicts the fact that $y <_m z$.

\begin{figure}
\begin{center}
% https://tikzcd.yichuanshen.de/#N4Igdg9gJgpgziAXAbVABwnAlgFyxMJZABgBoBGAXVJADcBDAGwFcYkQAPEAX1PU1z5CKcqWLU6TVuwBePPiAzY8BImXE0GLNohABPef2VCiAZgoSt03V15HBqlOY2Tt7A3cUCVw5ABYxSykdEDlPJQdfMj8gtxtDL2NHElJTWOt9BIifIlE0zWDZLO8TJ1IYgrjOYqTfc3zXDI8FbNL-VPSQsIkYKABzeCJQADMAJwgAWyQAJhocCCRRRpCpzzHJpDIQecXKjIBrEBpGegAjGEYABRLHEFGsPoALHAT1qcRZ7YXELasQxgAOgDGDAAI4AAkOxzOF2utXY9yeLzW43enx2iCW5zAUCQADYtowsGAQlB6HBHr0jst2IxqSdzlcbsIQMTsLBXqjdl9uYwIBA0EQAJxkYZMOAwCQM2HM9hsrAc47E0nkym4vb-TkbRDmHmIACsGvYUJA0qZ8N0iOeWveAT1hppulWCjeSDtGN1f1p4IAfJD6TDzZEEQ9rSjte7vg7sbjEASlST2GSKVSjboTWa4cHdPKOeH3g6MQ6+QLhaLxZLoYyszkc2B2WwEyqU+rHSBDvn8XNvgAONPtgPV2WW0PIl1cxBC7uLJZe9ODmUWu6jm2LLYY8izwpO1eIADs08QfbbzpGE7xh4PbcY4IAPP7O3HD8eY0gALR+Kem5VJ1Wp68LkGtasvWCpsI+x4YseJaCigIqkGKjASlKgY1qUIENvSP66MmarUnOpq7l+G6-Nupq+g+47auQ67fJuNCvnGpFVBm2EgLh-6ZsOGFgbu9F6vxMFlghFYoUOS65o236Jjhf6tgRHaUNwQA
\begin{tikzcd}
y \arrow[r, "m"] \arrow[d, "l\leq k"'] \arrow[r, "l" description, dashed, bend left=60]  & z \arrow["l" description, dashed, loop, distance=2em, in=35, out=325] &  & y \arrow[r, "m"] \arrow[d, "l > k"'] \arrow[r, "k" description, dashed, bend left=60] & z \arrow["k" description, dashed, loop, distance=2em, in=35, out=325] \\
x \arrow[ru, "k"']                                                                       &                                                                       &  & x \arrow[ru, "k"']                                                                    &                                                                       \\
                                                                                         &                                                                       &  &                                                                                       &                                                                       \\
y \arrow[r, "m"]                                                                         & z \arrow["l" description, dashed, loop, distance=2em, in=35, out=325] &  & y \arrow[r, "m"] \arrow[r, "k" description, dashed, bend left=60]                     & z \arrow["k" description, dashed, loop, distance=2em, in=35, out=325] \\
x \arrow[ru, "k"'] \arrow[u, "l < k"] \arrow[ru, "l" description, dashed, bend right=49] &                                                                       &  & x \arrow[ru, "k"'] \arrow[u, "l > k"]                                                 &                                                                      
\end{tikzcd}
\end{center}
\caption{Impossible relations in the proof of Lemma \ref{LemmaStratified}.}\label{FigureStrat}
\end{figure}
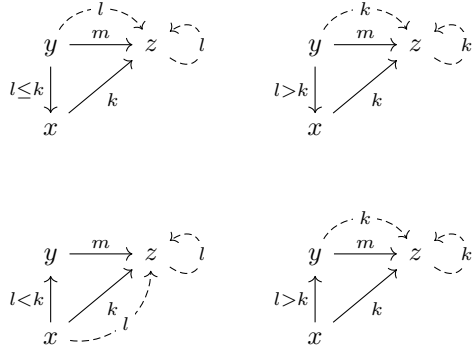

Thus, we have $x <_l y$. If $l < k$, then we have $x <_l z$ by Lemma \ref{LemmaJFrame},  contradicting the fact that $x <_k z$. If $k < l$, then we have $y <_k z$ by Lemma \ref{LemmaJFrame}, contradicting the fact that $y <_m z$. These cases which have been ruled out are illustrated in Figure \ref{FigureStrat}.
Therefore, the only possibility is that $x <_k y$, which is what is needed for stratification.
\endproof

Stratified, tree-like $J$-frames are sometimes called \textit{$J$-trees}, and so one might refer to hereditarily linear $J$-frames as \textit{$J$-lines}.

Let $(T, <_0, \hdots, <_{n-1})$ be a structure with $n$ binary relations and $x,y\in T$.
We write $\equiv_{n}$ for identity on $T$ and, inductively, $\equiv_k$ for the reflexive and symmetric closure of the union of $<_k$ and $\equiv_{k+1}$
%\[x<_k z\leftrightarrow y<_k z.\]
We refer to equivalence classes under $\equiv_k$ as \textit{$k$-planes}. \label{DefKPlanes}

\begin{lemma} \label{LemmaHLChar}
Let $(T, <_0, \hdots, <_{n-1})$ be a structure with $n$ binary relations. Then,
the following are equivalent:
\begin{enumerate}
\item\label{LemmaHLChar1} $T$ is a hereditarily linear $J_n$-frame; and 
\item\label{LemmaHLChar2} $T$ forms a single $\equiv_0$ equivalence class and for each $k < n-1$, $\equiv_{k+1}$ is compatible with $<_k$ and $<_k$ is a conversly well-founded linear order on $(k+1)$-planes when restricted to any given $k$-plane. 

\end{enumerate}
\end{lemma} 
\proof
Suppose $T$ is hereditarily linear $J_n$-frame. Using the fact that $T$ is stratified we show that  $<_k$  compatible with $\equiv_{k+1}$, for each $k<n$. The order between $(k+1)$-planes within a $k$-plane is total since if $x$ and $y$ from the same $k$-plane belong to different $(k+1)$-planes, then by hereditary linearity we have either $x <_k y$ or $y <_k x$.

Conversely, first observe that the structure is a $J_n$-frame. Suppose $k < m < n$. By hypothesis, $x \equiv_{m} x'$ and $y \equiv_{m} y'$ imply that $x <_k x'$ if and only if $y <_k y'$. Thus, if $x <_k y <_m z$, then we have $y \equiv_m z$ and thus $x <_k z$. The other two conditions are verified similarly. 
To see that it is hereditarily linear, let $x,y \in T$ be different and let $k$ be least so that $x \not\equiv_{k+1} y$. Then, $x \equiv_{k} y$. By \eqref{LemmaHLChar2} we have one of $x <_k y$, $y <_k x$, or else $x \equiv_{k+1} y$. By induction on $k < n$, we find some $l < n$ such that $x <_l y$ or $y <_l x$.
\endproof

\begin{figure}[h]
\begin{center}
\includegraphics[scale=0.08]{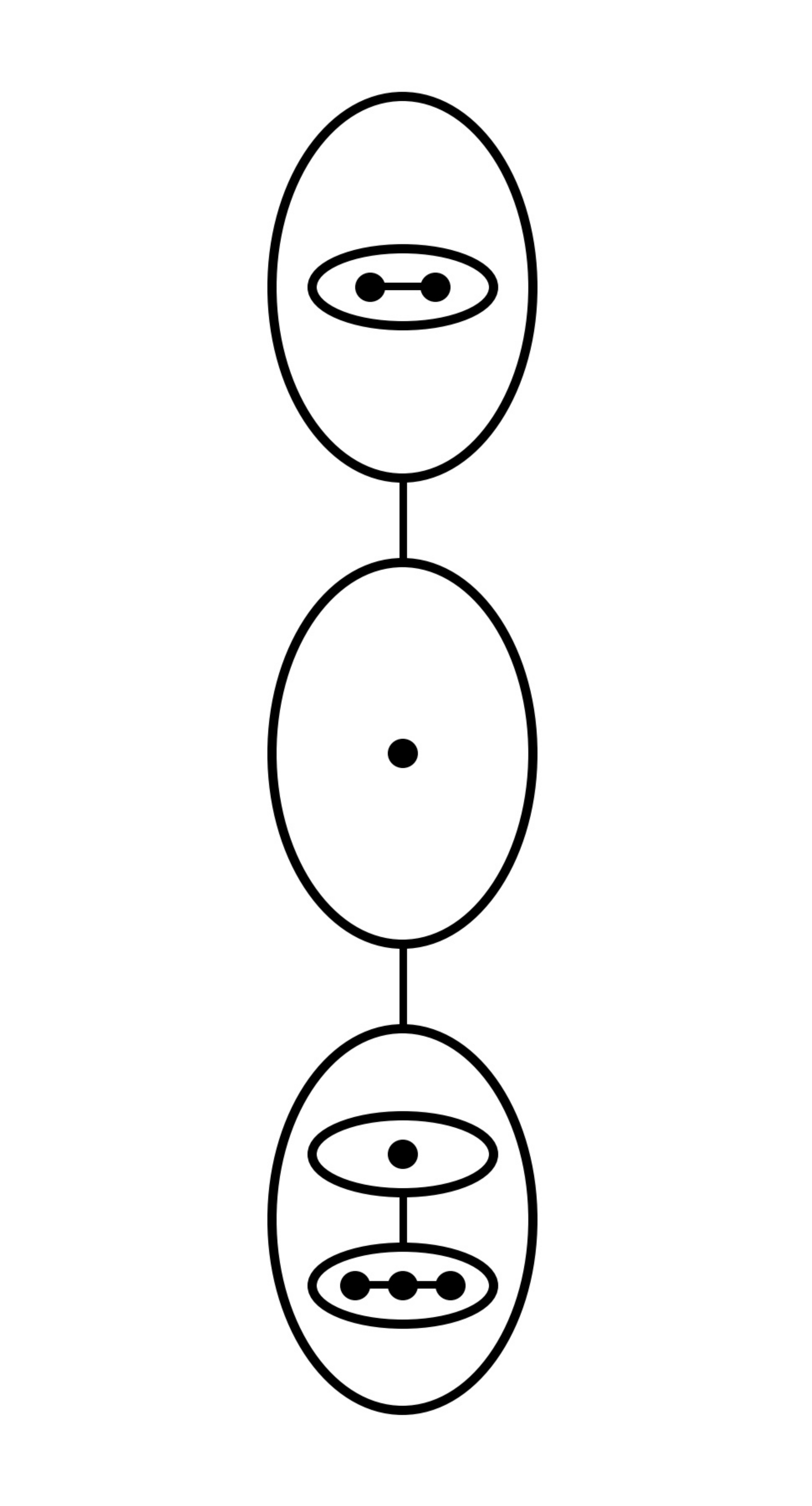}
\end{center}
\caption{A hereditary linear $J_2$-frame.}\label{FigureHL}
\end{figure}

See Figure \ref{FigureHL} for an illustration of Lemma \ref{LemmaHLChar}.

\subsection{Hereditarily Linear Semantics}
\begin{theorem}\label{TheoremCompletenessJlin}
$\ja_n.\lin$ is sound and complete for the class of finite, hereditarily linear  $J_n$-frames.
\end{theorem}
\proof
Soundness is easy to prove similarly to Theorem \ref{TheoremValidity}. Nonetheless, we include a direct proof for convenience, and then move on to completeness.\smallskip

\noindent\textsc{Proof of Soundness.} It is shown in Beklemishev \cite{Be04} that all $J_n$-frames satisfy all axioms of $\ja_n$, so we only need to check that the pseudo-linearity axiom is valid. 
\begin{lemma}\label{LemmaSoundness}
Every instance of pseudo-linearity is valid in every hereditarily linear, stratified $J_n$-frame.
\end{lemma}
\proof
We fix $m<k$ and a world $w\in T$ and verify that $m$ pseudo-linearity
\begin{align*}
\Vdash [m] \bigg(\Big(\bigwedge_{m\leq k< n}[k] \phi \Big) \to \psi \bigg) 
\vee 
[m] \bigg(\Big(\bigwedge_{m\leq k< n}[k]\psi \wedge \psi \Big) \to \phi \bigg) 
\end{align*}
holds in $w$ for any $\phi$ and $\psi$ and any interpretation of propositional variables. For this it is enough  to verify that it holds in $w$ as a world of the subframe $C$  that consists of the worlds reachable from $w$, first going by $<_m$ and then by relations $<_k$, for $m\leq k$; since $T$ is a $J_n$-frame, $C$, of course, is simply the $<_m$-cone of $w$. In $C$ points are $<_m$-reachable from $w$ if and only if they are $\bigcup_{m\le k<n} <_k$ reachable. That is the validity of $m$ pseudo-linearity in $w$ as a world in $C$ is equivalent to the validity of
$$\Box (\Box \phi\to \psi)\lor \Box((\Box \psi\land \psi)\to \phi),$$
where $\Box$ is the modality corresponding to the relation $\bigcup_{m\le k<n} <_k$. But the latter formula is just axiom 3 which is well-known to hold in any linear frame. And using the equivalent characterization of the hereditary linear frames from Lemma \ref{LemmaHLChar} it is easy to see that $\bigcup_{m\le k<n} <_k$ is linear in $C$.
\endproof

\noindent\textsc{Proof of Completeness}.
Let $\phi$ be a formula consistent with $\ja_n.\lin$.
Throughout this proof, we consider only $\mathcal{L}_{[<n]}$-formulas  containing only propositional variables in $\phi$. The proof is a modified ``canonical model'' construction, though it is significantly more involved than is typical. The main reason is that a $J$-line $W$ must be both hereditarily linear and converse wellfounded and irreflexive. The first of these conditions seems to require that many of the elements of $W$ are related via the accessibility relations, while the other two conditions seem to require the opposite. The construction must implicitly deal with this tension.
The simplest way to carry out the proof seems to be to first consider the canonical model $W$. There is no reason to suppose that $W$ will thus be converse wellfounded or hereditarily linear, and in fact it will be neither. 
However, we shall be able to isolate a particular submodel $W^*$ of $W$ which is both; this will consist of particularly nice worlds $w \in W$ which we call \textit{$k$-maximal}. Since $W^*$ is a proper subframe of $W$, a further argument will then be required to ensure that $\phi$ is still satisfied in $W^*$, and this will finish the proof.

Thus, let $W$ be the set of all maximal consistent set of formulas.
For such sets $x,y$ and $m< n$, we put $x <_m y$ if the following two conditions hold:
\begin{enumerate}
\item\label{Cond1CanPre} Suppose that $[m]\psi$ is such that $[m]\psi \in x$. Then, $\psi \in y$.
\item\label{Cond3CanPre} 
There is $\psi$ such that $[ m]\psi \in y\setminus x$.
\end{enumerate}

\begin{lemma}\label{LemmaCond1bBekl}
Suppose that $m\leq k< n$ and let $x,y \in W$ be such that $x <_m y$ and $[m]\psi \in x$. Then, $[k]\psi \in y$.
\end{lemma}
\proof
Since 
$\ja \vdash [m]\psi \to [m][k]\psi$
and $x$ is maximal consistent, we must have $[m][k]\psi \in x$ and thus $[k]\psi \in y$ by condition \eqref{Cond1CanPre} of the definition of $x <_m y$.
\endproof

\begin{lemma}\label{LemmaCond2Bekl}
Suppose that $k < m< n$ and let $x,y \in W$ be such that $x <_m y$. Then, for each $\psi$ we have $[k]\psi \in x\leftrightarrow [k]\psi \in y$.
\end{lemma}
\proof
Suppose $[k]\psi \in x$. 
Since 
$\ja \vdash [k]\psi \to [m][k]\psi$
and $x$ is maximal consistent, we must have $[m][k]\psi \in x$ and thus $[k]\psi \in y$ by condition \eqref{Cond1CanPre} of the definition of $x <_m y$. Similarly, if $[k]\psi \not \in x$, then we have $\langle k\rangle \lnot \psi \in x$, as it is maximal consistent. Since 
$\ja \vdash \langle k\rangle \lnot \psi \to [m]\langle k\rangle \lnot \psi$,
we have $[m]\langle k\rangle \lnot \psi \in x$ and thus $\langle k\rangle \lnot \psi \in y$ as above; thus $[k]\psi \not \in y$ as desired.
\endproof

Beklemishev \cite{Be10} proved the relational completeness of $\ja$ for finite $J$-frames. In his proof, he constructs finite models for a formula consistent with $\ja$ via maximal consistent subsets of a given finite set of formulas. Thus, the proofs of Lemma \ref{LemmaCond1bBekl} and Lemma \ref{LemmaCond2Bekl} do not go through for his definition of the accessibility relations $<_m$ and so he incorporates them into their definition. In our case, they follow from our definition of $<_m$, as we have seen. Using them, he shows:

\begin{lemma}[Beklemishev, Lemma 1.1 of \cite{Be10}]  \label{LemmaJFrame}
$(W, \vec <)$ is irreflexive and satisfies the following conditions:
\begin{enumerate}
\item \label{LemmaJFrame1}$x <_k y$ and $x <_m z$ together imply $y <_m z$ whenever $m < k$;
\item \label{LemmaJFrame2}$x <_k y <_m z$ implies $x <_m z$ whenever $m \leq k$;
\item \label{LemmaJFrame3}$x <_m y <_k z$ implies $x <_m z$ whenever $m \leq k$.
\end{enumerate}
\end{lemma}
Lemma \ref{LemmaJFrame} asserts that $W$ is a $J_n$-frame, except possibly for the condition of converse wellfoundedness. 

We define the natural valuation on $W$ by putting 
\begin{equation}\label{eqValuationCano}
W,x\Vdash \psi \leftrightarrow \psi \in x
\end{equation}
whenever $\psi = p$ is a propositional variable. Once more via Lemma \ref{LemmaCond1bBekl} and Lemma \ref{LemmaCond2Bekl} we have access to Beklemishev's argument, which yields the following:

\begin{lemma}[Beklemishev, Lemma 1.2 of \cite{Be10}]\label{LemmaValuationCano}
The equivalence \eqref{eqValuationCano} holds for arbitrary formulas $\psi$.
\end{lemma}

The next step is to restrict $W$ to a finite, hereditarily linear submodel $W^*$ while preserving the satisfiability of formulas. The key is the following notion:
Given $k<n$, let us say that a set $x$ of formulas is \textit{$k$-maximal} if there is a formula $\psi$ such that 
\[\lnot\psi \wedge\bigwedge_{k\leq m< n}[m]\psi \in x.\]
Clearly, if $x$ is $k$-maximal, then it is $m$-maximal for all $k\leq m< n$.

\begin{lemma}\label{ZeroMaximalxstar}
Suppose that $\phi$ is consistent with $\ja_n.\lin$. Then, there is a $0$-maximal $x^* \in W$ such that $\phi \in x^*$.
\end{lemma}
\proof
We will find such a $x^*$ where $0$-maximality is witnessed by $\lnot\phi$.
By hypothesis, $\phi$ is consistent with $\ja_n.\lin$.  Let $x \in W$ be such that $\phi \in x$. If $[0]\lnot\phi \not\in x$, then we have 
$x \Vdash \lnot [0]\lnot\phi$
by Lemma \ref{LemmaValuationCano}.
By $\mathsf{GL}$ for $[0]$, we have 
\begin{equation}
\mathsf{GL}\vdash \lnot [0]\lnot \phi \to \langle 0\rangle (\phi\wedge [0]\lnot\phi).
\end{equation}
Thus, we can find $y$ with $x <_0 y$ and 
$y\Vdash \phi\wedge [0]\lnot\phi$.
By Lemma \ref{LemmaValuationCano}, we then have $[0] \lnot\phi \in y$.
If $[1] \lnot\phi\not\in y$, we can similarly find $z$ with $y <_1 z$ and 
$z\Vdash \phi\wedge [1]\lnot\phi$.
Thus, we have $\phi \wedge [1]\lnot\phi \in z$. By the definition of $<_1$, we also have $[0]\lnot\phi \in z$. Continuing this way, we eventually find some $w$ containing each of
$\phi, [0]\lnot\phi, \dots, [n-1]\lnot\phi$,
in which case we have 
\[\phi\wedge\bigwedge_{k< n} [k]\lnot\phi \in w,\]
by maximality.
\endproof

Similarly, and using the axioms of $\ja$, we obtain:
\begin{lemma}\label{LemmaConstructionWStar}
Suppose $x \in W$, $\psi \in \Sigma$, and $\lnot[k]\psi \in x$, then there is $y\in W$ such that $x<_k y$ and
\[\lnot\psi \wedge \bigwedge_{k\leq m< n}[m]\psi \in y.\]
In other words, $y$ is $k$-maximal, as witnessed by $
\psi$.
\end{lemma}
\iffalse
\proof
This is similar to Lemma \ref{ZeroMaximalxstar}, but we additionally need to verify that $x <_k y$.
By Lemma \ref{LemmaValuationCano} we can find $y$ with $x <_k y$ and $\lnot\psi \in y$. Since $y$ satisfies $\mathsf{GL}$ for all boxes $[k]$ and
\[\mathsf{GL} \vdash \lnot [k]\psi \to \langle k\rangle ([k]\psi \wedge \lnot\psi),\]
it follows that if $[k]\psi \not\in y$, then we can find $z$ such that $y <_k z$ and 
\[\lnot\psi \wedge [k]\psi \in z.\]
By transitivity (Lemma \ref{LemmaJFrame}) it follows that $x <_k z$ and thus we may assume without loss of generality that 
\[\lnot\psi \wedge [k]\psi \in y.\]
Inductively, if $[m]\psi \not \in y$, then we have
\[\mathsf{GL} \vdash \lnot [m]\psi \to \langle m\rangle ([m]\psi \wedge \lnot\psi),\]
so we can find $z$ with $y<_m z$ and 
\[\lnot\psi \wedge [m]\psi \in z.\]
Since $k <m$ and applying Lemma \ref{LemmaJFrame}, we have 
\[(m' < m \wedge [m']\psi \in y)\to [m']\psi \in z.\]
Moreover, we also have $x<_k z$ by Lemma \ref{LemmaJFrame}. Continuing this way, we find some $z$ such that $x<_k z$ and 
\[\lnot\psi \wedge \bigwedge_{k\leq m\leq n}[m]\psi \in z,\]
as desired.
\endproof 
\fi

The following is the key lemma in the construction:
\begin{lemma}\label{LemmaPseudoLin}
Suppose that $x$, $y$, $z$, and $k < n$ are such that $z<_k x$, $z <_k y$ and both $x$ and $y$ are $k$-maximal. Then, one of the following holds:
\begin{enumerate}
\item $x= y$,
\item $x <_m y$ for some $k\leq m <n$, or 
\item $y <_m x$ for some $k\leq m <n$.
\end{enumerate}
\end{lemma}
\proof
We will use Lemma \ref{LemmaValuationCano} without mention.

Suppose that $x,y,z$ form a counterexample to the statement of the theorem. By the fact that $x \neq y$, we obtain some formula $\chi$ such that $\chi \in y\setminus x$, so that:
\[y\Vdash \chi; \quad x\Vdash \lnot\chi.\]

\begin{claim}
For each $k\leq m< n$, there is some $\psi^x_m$ such that $[m]\psi^x_m \in x$ and $\psi^x_m \not\in y$. Similarly,  there is some $\psi^y_m$ such that $[m]\psi^y_m \in y$ and $\psi^y_m \not\in x$.
\end{claim}
\proof
By hypothesis, $y$ is $k$-maximal, so there is a formula $\psi^x_m = \psi$ such that 
\begin{equation}\label{eqYmaximalClaimLin}
\lnot\psi \wedge \bigwedge_{k \leq l < n} [l] \psi \in y.
\end{equation}
In particular, we have $[m]\psi \in y$ and thus $[m]\psi \in x$. However, by \eqref{eqYmaximalClaimLin} we have $\lnot\psi \in y$, so $\psi \not\in y$, as desired. The formula $\psi^y_m$ is obtained similarly.
\endproof

Define formulas $\psi^x$ and $\psi^y$ by 
\begin{align*}
\psi^x = \bigvee_{k\leq m< n} \psi^x_m, \quad
\psi^y = \bigvee_{k\leq m< n} \psi^y_m \vee \chi.
\end{align*}
\begin{claim}\label{ClaimQLinPsixy}
The following hold:
\begin{enumerate}
\item $\psi^x \not\in y$,
\item $[m]\psi^x \in x$ for each $m$ with $k\leq m < n$,
\item $\psi^y \not\in x$,
\item $[m]\psi^y \in y$ for each $m$ with $k\leq m < n$,
\item $\psi^y \in y$.
\end{enumerate}
\end{claim}
\proof
Since $\psi^x_m\not\in y$ for all $m$, we also have $\psi^x \not\in y$. Similarly, $\psi^y_m\not\in x$ and $\chi \not\in x$, so $\psi^y \not\in x$. For each $m$, we have $[m]\psi^x_m \in x$ and thus $[m]\psi^x \in x$. Similarly, we have $[m]\psi^y \in y$. Finally, $\chi \in y$, so $\psi^y\in y$.
\endproof

Now, since $z$ is consistent with $\ja_n.\lin$ and maximal, it contains all instances of the quasi-linearity axiom. In particular:
\begin{align}
z\Vdash 
[k] &\bigg( \Big(\bigwedge_{k\leq m< n} [m] \psi^x\Big) \to \psi^y\bigg) \vee [k] \bigg( \Big( \bigwedge_{k\leq m< n} [m]\psi^y \wedge \psi^y \Big) \to \psi^x\bigg).
\label{eqZBigDisj}
\end{align}
There are two cases. Suppose first that $z$ satisfies the first disjunct in \eqref{eqZBigDisj}. Since $z <_k x$, we have 
\begin{align*}
x\Vdash \Big(\bigwedge_{k\leq m< n} [m] \psi^x\Big) \to \psi^y.
\end{align*}
However by Claim \ref{ClaimQLinPsixy} we have $x\not\Vdash \psi^y$ while $x \Vdash [m]\psi^x$ for all $m$ with $k\leq m< n$, so we have a contradiction. 

The second case is that $z$ satisfies the second disjunct in \eqref{eqZBigDisj}. Since $z <_k y$, we have 
\begin{align*}
y\Vdash \Big( \bigwedge_{k\leq m< n} [m]\psi^y \wedge \psi^y \Big) \to \psi^x
\end{align*}
However by Claim \ref{ClaimQLinPsixy} we have $y\not\Vdash \psi^x$ while $y \Vdash [m]\psi^y$ for all $m$ with $k\leq m< n$ and moreover $y\Vdash \psi^y$, so we have a contradiction once more.  This completes the proof of the lemma.
\endproof

By hypothesis $\phi$ is consistent with $\ja_n.\lin$.
Let $\Sigma$ be the smallest set of formulas such that
\begin{enumerate}
\item $\Sigma$ contains all subformulas of $\phi$,
\item if $\psi \in\Sigma$ and the outermost connective of $\psi$ is not a negation, then, $\lnot\psi \in\Sigma$,
\item if $[k]\psi\in\Sigma$ for some $k$ and $k'< n$, then $[k']\psi\in\Sigma$.
\end{enumerate}
 
Let $x^*$ be a $0$-maximal element of $W$ such that $\phi\in x^*$. This exists by Lemma \ref{ZeroMaximalxstar}.
%Below, we denote by $\Sigma$ the collection of subformulas of \eqref{eqBigWedge}. 
Write $W^*_{-1} = \{x^*\}$.
We let $W^*_0$ consist of $x^*$ and of all $y$ such that $x^*<_0 y$ and which are $0$-maximal as witnessed by some $\psi =: \psi^y$ such that $\lnot [0]\psi \in x^* \cap \Sigma$.

Inductively, let $W^*_{k}$ consist of all elements $x\in W^*_{k-1}$ and of all $y\in W$ such that $x <_k y$ for some $x \in W^*_{k-1}$ and $y$ is $k$-maximal as witnessed by some $\psi =: \psi^y$ such that $\lnot [k]\psi \in x \cap \Sigma$.

Finally, we write 
\[W^* = \bigcup_{k< n} W^*_k.\]
We let the \textit{degree} of $y \in W^*$ be the least $k$ such that $y \in W^*_k$. We regard $(W^*, \vec <)$ as a subframe of $W$ induced by the relations $\vec <$.

\begin{lemma}\label{LemmaHereditaryLinear}
$(W^*, \vec <)$ is hereditarily linear.
\end{lemma}
\proof
This is proved by induction. $W^*_{-1}$ is trivially hereditarily linear. Suppose $W^*_{k-1}$ is and let $u,v \in W^*_{k}$. Then, \smallskip

\noindent\textsc{Case I}: $u,v \in W^*_{k-1}$. Immediate.\smallskip

\noindent\textsc{Case II}: $u \in W^*_{k-1}$, $v \not \in W^*_{k-1}$. Then, $v$ is $k$-maximal.
Let $x \in W^*_{k-1}$ be such that $x <_k v$ witnesses the fact that $v \in W^*_{k}$. If $x = u$, then we have $u <_k v$ by hypothesis.
Otherwise inductively we have  either $u <_m x$ for some $m$, in which case $u <_m v$ by Lemma \ref{LemmaJFrame} if $m < k$, $u <_k v$ by Lemma \ref{LemmaJFrame} if $k \leq m$;
or else $x <_m u$ for some $m$, in which case similarly we have $v <_m u$ by Lemma \ref{LemmaJFrame} if $m < k$, $u <_k v$ by Lemma \ref{LemmaJFrame} if $k < m$, or $u$ and $v$ are $l$-comparable for some $l\geq k$ by Lemma \ref{LemmaPseudoLin} if $k = m$. This application of Lemma \ref{LemmaPseudoLin} is justified by the fact that $v$ is $k$-maximal and $u$ is $k'$-maximal for some $k' < k$ and thus also $k$-maximal.
\smallskip

\begin{figure}[h]
\begin{center}
% https://q.uiver.app/#q=WzAsOCxbMywxLCJ4Il0sWzQsMSwieSJdLFszLDAsInUiXSxbNCwwLCJ2Il0sWzAsMSwieCJdLFsxLDEsInkiXSxbMCwwLCJ1Il0sWzEsMCwidiJdLFswLDEsIm0iLDJdLFswLDIsImsiXSxbMSwzLCJrIiwyXSxbNCw1LCJtIiwyXSxbNCw2LCJrIl0sWzUsNywiayIsMl0sWzQsNywibSIsMix7InN0eWxlIjp7ImJvZHkiOnsibmFtZSI6ImRhc2hlZCJ9fX1dLFs2LDcsIm0iLDAseyJzdHlsZSI6eyJib2R5Ijp7Im5hbWUiOiJkYXNoZWQifX19XSxbMywyLCJtJyAoXFxnZXEgaykiLDIseyJjdXJ2ZSI6Mywic3R5bGUiOnsiYm9keSI6eyJuYW1lIjoiZGFzaGVkIn0sImhlYWQiOnsibmFtZSI6Im5vbmUifX19XSxbMCwzLCJrIiwyLHsic3R5bGUiOnsiYm9keSI6eyJuYW1lIjoiZGFzaGVkIn19fV1d
\[\begin{tikzcd}
	u & v && u & v \\
	x & y && x & y
	\arrow["m"', from=2-4, to=2-5]
	\arrow["k", from=2-4, to=1-4]
	\arrow["k"', from=2-5, to=1-5]
	\arrow["m"', from=2-1, to=2-2]
	\arrow["k", from=2-1, to=1-1]
	\arrow["k"', from=2-2, to=1-2]
	\arrow["m"', dashed, from=2-1, to=1-2]
	\arrow["m", dashed, from=1-1, to=1-2]
	\arrow["{m' (\geq k)}"', curve={height=18pt}, dashed, no head, from=1-5, to=1-4]
	\arrow["k"', dashed, from=2-4, to=1-5]
\end{tikzcd}\]
\end{center}
\caption{Comparability of $u$ and $v$ in the proof of Lemma \ref{LemmaHereditaryLinear}, \textsc{Case III}, assuming $m < k$ (left) and $k \leq m$ (right).}
\label{FigureHLProof}
\end{figure}
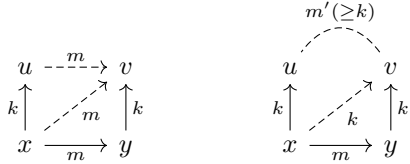

\noindent\textsc{Case III}: $u,v \not\in W^*_{k-1}$. Then, $u$ and $v$ are $k$-maximal. Let $x <_k u$ and $y <_k v$ witness the facts that $u,v \in W^*_{k}$. If $x = y$, then Lemma \ref{LemmaPseudoLin} yields the $m$-comparability of $u$ and $v$ for some $m \geq k$. Otherwise inductively $x <_m y$ or $y <_m x$ for some $m$. Let us consider the first case; the second one is analogous. This situation is illustrated in Figure \ref{FigureHLProof}.

If $m < k$, then $x <_m v$ by  Lemma \ref{LemmaJFrame} and thus $u <_m v$ by Lemma \ref{LemmaJFrame}. If $k \leq m$, then $x <_k v$ by Lemma \ref{LemmaJFrame}, and hence $u$ and $v$ are $m'$-comparable for some $m'\geq k$ by Lemma \ref{LemmaPseudoLin}. Note that this application of Lemma \ref{LemmaPseudoLin} is justified by the fact that $u$ and $v$ are both $k$-maximal.
\endproof

%Thus, $x^*$ is the least element of $W^*$ in the sense of all accessibility relations $<_k$.

The following lemma asserts that $y\in W^*$ as above is uniquely determined by the formula $\psi^y$, locally.
\begin{lemma}\label{LemmaFiniteness}
Let $y$ and $z$ be elements of $W^*$ of degree $k$ and let $x\in W^*_{k-1}$ be such that $x <_k y,z$. If $\psi^y = \psi^z$, then $y = z$.
\end{lemma}
\proof
Since $x <_k y,z$, Lemma \ref{LemmaPseudoLin} gives three possibilities:
The first is that for some $m$ with $k\leq m$ we have $y <_m z$. However, $y \Vdash [m]\psi^y$ and $z \Vdash \lnot\psi^z$, so this is impossible. The second is that  for some $m$ with $k\leq m$ we have $z <_m y$. This is impossible for the same reason. The third one is that $y = z$ as desired.
\endproof

\begin{lemma}\label{LemmaWStarFinite}
$W^*$ is finite.
\end{lemma}
\proof
Immediate from Lemma \ref{LemmaFiniteness}.
\endproof

$W^*$ is a subframe of $W$. 
We now also regard it as a submodel of $W$ with the same interpretations for all variables. Henceforth we use $y\Vdash \psi$ to mean $W^*, y\Vdash \psi$, though we may clarify if necessary. We need to show that the conclusion of Lemma \ref{LemmaValuationCano} remains true for $W^*$.

\begin{lemma}\label{LemmaValuationCanoStar}
For all formulas $\psi \in \Sigma$ and all $y \in W^*$, we have 
\begin{equation*}
W^*,y\Vdash \psi \leftrightarrow \psi \in y.
\end{equation*}
\end{lemma}
\proof
By induction on $\psi$. The crucial case is that of formulas of the form $[m]\psi$. If $y\not\Vdash[m]\psi$, so that $z\not\Vdash \psi$ for some $z\in W^*$ with $y <_m z$, then in particular we have 
$W,y\Vdash \lnot[m]\psi$,
and hence $\lnot[m] \psi \in y$, by Lemma \ref{LemmaValuationCano}, so $[m]\psi\not\in y$. For the converse, suppose $[m]\psi\not\in y$. We distinguish two cases. Let $k$ be the degree of $y$.\smallskip

\noindent\textsc{Case I}: $k < m$. We have $y \in W^*_k$. By Lemma \ref{LemmaConstructionWStar}, there is $z \in W^*_m$ such that $y<_m z$ and $\lnot\psi \in z$. By induction hypothesis, $z \Vdash \lnot\psi$, so $y\Vdash \lnot [m]\psi$.\smallskip

\noindent\textsc{Case II}: $m\leq k$.
Again we have $y\in W^*_k$. Thus, there is some $k_0 < k$ and $x_0 \in W^*_{k_0}$ such that $x_0 <_k y$. If $m\leq k_0$, find some $k_1 < k_0$ and some $x_1 \in W^*_{k_1}$ such that $x_1<_{k_0} x_0$. Continue this way, producing a sequence
\[k_l < k_{l-1} < \hdots < k_0 < k,\]
where $l$ is minimal such that $k_l < m$, and worlds $x_0, \hdots, x_l$ such that $x_i \in W^*_{k_i}$ and $x_{i+1} <_{k_i} x_i$ for each $i<l$.
We must have $\lnot [m] \psi \in x_l$, for otherwise $[m]\psi \in x_l$ by maximality and thus inductively we would be able to show that $[m]\psi \in y$, using Lemma \ref{LemmaCond2Bekl} and the fact that 
\[m \leq k_{l-1} < k_{l-2} < \dots k_0 < k,\] 
contradicting the fact that $[m]\psi\not \in y$. Thus, $\lnot [m] \psi \in x_l$. 
Applying Lemma \ref{LemmaJFrame} repeatedly, we inductively see that $x_l <_{k_{l-1}} x_{l-1-i}$ for each $0\leq i < l$ and then that $x_l <_{k_{l-1}} y$. 

Since $x_l \in W^*_{k_l}$ and $k_l < m$, Lemma \ref{LemmaConstructionWStar} provides some $z\in W^*_m$ such that $x_l <_m z$ and 
\[\lnot\psi \wedge \bigvee_{m\leq l < n}[l]\psi \in z.\] 
By induction hypothesis, this implies in particular that $z \Vdash \lnot\psi$.
We must now distinguish two cases:\smallskip

\noindent\textsc{Case IIa}: $m < k_{l-1}$. 
Since $x_l <_m z$, $x_l <_{k_{l-1}} y$, and $m < k_{l-1}$,  Lemma \ref{LemmaJFrame} yields $y <_m z$. It follows that $y\Vdash \lnot[m]\psi$, as desired.\smallskip

\noindent\textsc{Case IIb}: $m = k_{l-1}$. Thus, we have $x <_m z$ and $x <_m y$.
Applying Lemma \ref{LemmaPseudoLin} to $z$ and $y$, we obtain one of three possibilities: 

The first one is that $z = y$. However, we have $[m]\psi \in z$ and $\lnot[m]\psi \in y$, so this is impossible. The second is that $z <_l y$ for some $l$ with $m \leq l$. Since $[m]\psi \in z$ and $\lnot[m]\psi \in y$, we cannot have $m < l$, Lemma \ref{LemmaCond2Bekl}, so we must have $m = l$.
However, we have $[m]\psi \in z$, so  $z <_m y$ would imply $[m]\psi \in y$ by Lemma \ref{LemmaCond1bBekl}, which again is impossible. 

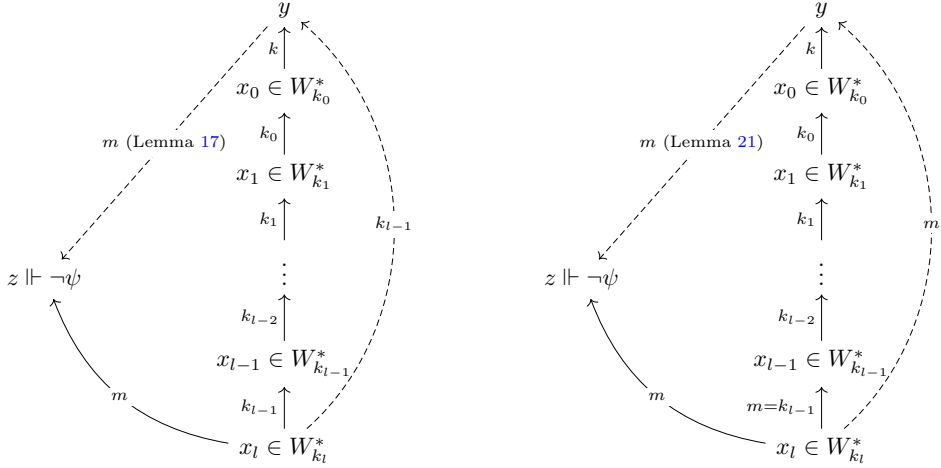
\begin{figure}[h]
\adjustbox{scale=0.9,center}{%
%\begin{center}
% https://tikzcd.yichuanshen.de/#N4Igdg9gJgpgziAXAbVABwnAlgFyxMJZAJgBoAWAXVJADcBDAGwFcYkQAPAfWEYFoAjAF8ABAB0xWMCIDqAPQBUPANY9+woSCGl0mXPkIoyAZmp0mrdhNpQIOBNt3Y8BImWJmGLNok5cB4pLS8krAqsJaOiAYzgZupAKeFj5+AAyBUrKKKlypmo7Rei6GJKSpSd7sAJ6RTvquRqQArBWWvtyMGcHZYVyM+VEx9SWppKY0Xm0gAF4SAGpQ9HAAFoGMkDgSaNi1hbENyADsza0pHV1Zoar9u0PFRMdUE8ns3LyCohKZITnvGrdFOIoY7jcyVXzWWz2AH7ErHDzPcF+AJfbpXfwDOr3YEJU6vXIXH69PIw4YPMp43w1Ap3IHIJpjSkzeaLFYSdZ2LY7IRmGBQADm8CIoAAZgAnCAAWyQoxAOAgSAEiKmqnexExIHFUsVNHlSDIYJV-l2WuliANesQxmVKVUqRNErN1rlCsQ5Bt7GUDu1iAZLplHt8qvUGtNSD9lvdIAARjAwFAkHxyABOGiMKQpVnLPkgQMgYMfXMgRj0WOMAAKgIaICk2Fg3rNEddADY89K06WYBWq4Ya2A62waLH40hjHkomG3bqW2mM+wszm24EcDAODhgCIABQAGRgksl9ECYpgIuAu-39AAUgAxMX0SUwIQASk0HbLldh7FrWHrBUnAiVf1EAEA1JltBtFVlS0ALzO0IMQY4gIADiXABeEQCwiP9HSQFCgNTQ1bTUPh1XggjoNlMDPWNbCfUQ6CDXTMBMyWbMEyHOME0QJMCKo3x22LTtu0-Xxv1-CccIQ6dFWdEt3x7L9+x-QcY040dKJefj4JA6TgNkudfAXdjCPYSVl1Xdct3PA8jxPM89wPcs4BgZhbG3KRn1fQT5JEvsBy0SghCAA
\begin{tikzcd}[scale=0.50]
                  &  & y \arrow[llddd, "m \text{ (Lemma \ref{LemmaJFrame})}" description, dashed]                                                                  &  &  &                  &  & y \arrow[llddd, "m \text{ (Lemma \ref{LemmaPseudoLin})}" description, dashed]                                                             \\
                  &  & x_0 \in W^*_{k_0} \arrow[u, "k"]                                                                                                            &  &  &                  &  & x_0 \in W^*_{k_0} \arrow[u, "k"]                                                                                                          \\
                  &  & x_1 \in W^*_{k_1} \arrow[u, "k_0"]                                                                                                          &  &  &                  &  & x_1 \in W^*_{k_1} \arrow[u, "k_0"]                                                                                                        \\
z\Vdash \lnot\psi &  & \vdots \arrow[u, "k_1"]                                                                                                                     &  &  & z\Vdash\lnot\psi &  & \vdots \arrow[u, "k_1"]                                                                                                                   \\
                  &  & x_{l-1} \in W^*_{k_{l-1}} \arrow[u, "k_{l-2}"]                                                                                              &  &  &                  &  & x_{l-1} \in W^*_{k_{l-1}} \arrow[u, "k_{l-2}"]                                                                                            \\
                  &  & x_l \in W^*_{k_l} \arrow[u, "k_{l-1}"] \arrow[uuuuu, "k_{l-1}" description, dashed, bend right=49] \arrow[lluu, "m" description, bend left] &  &  &                  &  & x_l \in W^*_{k_l} \arrow[u, "m = k_{l-1}"] \arrow[uuuuu, "m" description, dashed, bend right=49] \arrow[lluu, "m" description, bend left]
\end{tikzcd}
%\end{center}
}
\caption{Finding an $m$-extension $z$ of $y$ satisfying $\lnot\psi$ in the proof of Lemma \ref{LemmaValuationCanoStar} in cases \textsc{IIa} (left) and \textsc{IIb} (right).}\label{FigLemmaValuation} 
\end{figure}

The third possibility is that $y <_l z$ for some $l$ with $m\leq l$. As before, however, it follows from Lemma \ref{LemmaCond2Bekl} that $m = l$. If so, then $z$ witnesses the fact that $y\Vdash \lnot [m]\psi$, as desired.  See Figure \ref{FigLemmaValuation} for the picture. This proves the lemma.
\endproof

By Lemma \ref{LemmaWStarFinite}, $(W^*, \vec <)$ is finite and hence converse wellfounded.  By Lemma \ref{LemmaJFrame} $(W^*, \vec <)$ is thus a $J_n$-frame, and by Lemma \ref{LemmaHereditaryLinear} it is hereditarily linear. %By Lemma \ref{LemmaStratified}, $W^*$ is stratified. 
Finally, by Lemma \ref{LemmaValuationCanoStar} we have
$x^*\Vdash \phi$.
We have shown that $(W^*, \vec <)$ is a finite and hereditarily linear $J_n$-frame satisfying $\phi$. This proves the completeness theorem.
\endproof

It might be worth noting that the proof of the completeness theorem yields the following bound on the size of the smallest model of a formula.
\begin{corollary}
Suppose that $\phi$ is a formula which is consistent with $\ja_n.\lin$ and let $k$ be the number of boxes and diamonds occurring in $\phi$. Then, $\phi$ is satisfied in a hereditarily linear $J_n$-frame $(T, \vec <)$ such that
\[|T| \leq (k+1)\cdot k^{n-1}.\]
\end{corollary}
\proof
It suffices to inspect the construction of the model $W^*$ in the proof of the theorem. It begins with the root $x^*$ which will eventually satisfy $\phi$ and iteratively produces the models $W^*_l$. Each step adds for each world  $x \in W^*_l$ a neighbor for each formula $\lnot [l] \psi \in x\cap \Sigma$. 
For $\lnot [l] \psi$ to belong to $\Sigma$, it is necessary and sufficient that a formula sharing its negation-normal form with $[l']\psi$ be a subformula of $\phi$, for some $l' < n$. 
Since there are at most $k$ such formulas, and the index $l$ ranges over numbers $0 \leq l < n$, the result follows.
\endproof

\section{Pointwise Definability of Beklemishev-Icard Projections}\label{SectDefinability}
Ignatiev \cite{Ig93} introduced a universal frame for $\glp$ which we denote $\ig$. Beklemishev \cite{Be10} introduced a more general kind of model for $\glp$ via the so-called \textit{blow-up} construction. 
Beklemishev showed that if $T$ is a $J$-frame with hereditary root $r$ and $\phi$ is a $\mathcal{L}_{[<\omega]}$ formula such that $r, T\Vdash \phi\wedge M^+(\phi)$, where $M^+(\phi)$ is a certain conjunction of $\glp$-consequences of $\phi$ defined below, then there is a projection $\pi$ onto $T$ from its \textit{$\omega$-blow-up} $B(T)$ with the property that if $\pi(x) = r$, then $x, B(T) \Vdash \phi$.
Icard \cite{Ic08,Ic11} connected the two models by showing that when Beklemishev's blow-up method is applied to hereditarily linear $J$-frames, one obtains an initial segment of Ignatiev's frame. %Let us state the goal of this section before recalling the relevant definitions:

For each $\mathcal{L}_{[<n]}$-formula $\phi$  we define
\[M(\phi) = \bigwedge \bigg\{\bigwedge_{k = m + 1}^{n-1}\big([m]\psi \to [k]\psi\big) : [m]\psi \text{ is a subformula of } \phi \bigg\}\]
and moreover we write $M^+(\phi) = M(\phi) \wedge \bigwedge_{m < n}[m]M(\phi)$.

Ignatiev's frame $\ig$ is defined as the set of all sequences of ordinals $\{\alpha_i: i\in\mathbb{N}\}$ such that for all $i\in\mathbb{N}$ we have $\alpha_i < \varepsilon_0$ and $\alpha_{i+1} \leq \log(\alpha_i)$, where $\log(\alpha + \omega^\beta) := \beta$. %For technical reasons, we remove from $\ig$ the constant-zero sequence.
Every element of $\ig$ is eventually equal to zero and is thus identified with the unique finite sequence with the same support. The accessibility relations are defined by $\vec\alpha <_m \vec\beta$ if and only if $(\alpha_0, \hdots, \alpha_{m-1}) = (\beta_0, \hdots, \beta_{m-1})$ and $\alpha_m < \beta_m$. The set $\ig_\iota$ is defined as the set of all $\vec\alpha \in \ig$ whose first component is smaller than or equal to $\iota$ (so according to this definition, $\ig_\iota$ has a last level in the sense of $<_0$). Each level $\iota$ of Ignatiev's frame contains a distinguished point $\delta_\iota := \vec \alpha$ satisfying $\alpha_{i+1} = \log(\alpha_i)$ for all $i\in\mathbb{N}$. The \textit{main axis} of Ignatiev's frame consists of precisely these distinguished points; it is wellordered by $<_0$ and has length $\varepsilon_0$. Ignatiev \cite{Ig93} proved that every closed $\mathcal{L}_{[<\omega]}$-formula consistent with $\glp$ is satisfied at a point on the main axis. Conversely, every point $\delta_\iota$ on the main axis is closed-definable, in the sense that there is a closed $\mathcal{L}_{[<\omega]}$-formula $\phi^{\delta_\iota}$ such that  $\vec\gamma,\ig\Vdash \phi^{\delta_\iota}$ if and only if $\vec\gamma = \delta_\iota$. It follows that the set $\ig\setminus \ig_\iota$ can be characterized as consisting precisely of those points which satisfy $\langle 0 \rangle \phi^{\delta_\iota}$.

\begin{lemma}\label{LemmaDefinability}
  Suppose that $T$ is a finite, hereditarily linear $J$-frame. Then there is an ordinal $\iota$ and a projection $\pi: \ig_\iota \to T$ with the following properties:
  \begin{enumerate}
   \item \label{LemmaDefinability1} for any $t\in T$, the set $\pi^{-1}(t)$ is definable in $\ig$;
   \item \label{LemmaDefinability2} the root $\delta_\iota$ of $\ig_\iota$ is the only point mapped by $\pi$ to the root $r$ of $T$;
   \item \label{LemmaDefinability3} for any polymodal formula $\phi$, and any evaluation $\mathsf{ev}$ of its variables in $T$, if $M^+(\phi)$ holds at $r$ under this evaluation, then $\phi$ holds at $(\alpha_0,\ldots)\in \ig_\iota$ under the evaluation $\mathsf{ev}'$, where $\mathsf{ev}'$ is defined by $\mathsf{ev}'(p)=\pi^{-1}[\mathsf{ev}(p)]$ if and only if $\phi$ holds at $\pi((\alpha_0,\ldots))$ under the evaluation $\mathsf{ev}$;
   \item \label{LemmaDefinability4}as a consequence of (1) and (3), if $M^+(\phi)\land \phi$ is satisfiable at the root of $T$, then there is some closed substitution $\phi^*$ such that $\ig,\delta_\iota\Vdash \phi^*$.
  \end{enumerate}
\end{lemma}

 Our construction of $\pi$ -- while less general than Beklemishev's and tailored to our context -- is arguably simpler, and implicit in \cite{Be10} and in \cite{Ic08,Ic11}. %We then use it to prove that the conclusion of the lemma holds, by induction. 

\begin{remark}
The conclusion of the lemma also holds for the topological models of $\glp$ constructed by Beklemishev and Gabelaia \cite{BG13} and the $J$-map $j$ constructed in their proof instead of the projection $\pi$; this can be proved by induction on the construction of $j$ following the proof of \cite{BG13}. These models could be employed in place of Ignatiev's frame throughout our proof.
\end{remark}

\proof[Proof of Lemma \ref{LemmaDefinability}]
By simultaneous induction on $n$, given a hereditarily linear $J_n$-frame $T$, we construct an ordinal $\iota$ and projection $\pi: \ig_\iota \to T$ as well as veirfy the properties (\ref{LemmaDefinability1}) and (\ref{LemmaDefinability2}).
\smallskip

\noindent\textsc{Case I}: $T$ is a singleton. Then, $\iota = 0$ and the projection maps the unique point in $\ig_0$ to the unique point in $T$. The defining formula is simply $[0] \bot$. Validity of all the properties is straightforward to verify.\smallskip

\noindent\textsc{Case II}: $<_0^T$ is empty but $T$ is not a singleton. Letting $S = (T, <_1^T, \hdots, <_{n-1}^T)$, $S$ is a hereditarily linear $J_{n-1}$-frame. Let $\iota$, $\pi$, and defining formulas for $\pi^{-1}(t)$, $t\in T$, be given by the induction hypothesis applied to $S$. 
We define $\hat \iota$, $\hat \pi$, and defining formulas for $\hat \pi^{-1}(t)$, $t\in T$, for the case $T$. 
Let $\hat\iota=\omega^\iota$ and $\hat \pi ( (\alpha_0,\alpha_1,\ldots))=\pi((\alpha_1,\alpha_2,\ldots))$, and given a defining formula  $\psi$ for $\pi^{-1}(t)$ we construct a formula $\hat\psi'$ by replacing each $[k]$ and each $\langle k\rangle$ in $\psi$ by $[k+1]$ and $\langle k+1\rangle$, respectively.
Notice that $\hat \psi'$ is valid in $(\alpha_0,\alpha_1,\ldots)$ if and only if $\psi$ is valid in $(\alpha_1,\ldots)$. Hence $\hat\pi^{-1}(t)$ is defined by:
$$\hat\psi := \hat\psi'\land [0]\lnot \phi^{\hat\iota}.$$

Clearly, $\delta_\iota$ is the only point projected to the root of $T$.

\noindent\textsc{Case III}: $<_0^T$ is nonempty. Then, $T$ is a set of $1$-planes linearly ordered by $<_0^T$, say $P_0, \hdots, P_l$. By induction hypothesis  we have ordinals $\alpha_0, \hdots, \alpha_l$ and projections $\pi_i: \ig_{\alpha_i} \to P_i$. We let $\iota = \alpha_0 + 1+\alpha_1 + \dots + 1 +\alpha_l$. Then, $\ig_\iota$ can be regarded as the ordered sum of $\ig_{\alpha_0}, \hdots, \ig_{\alpha_l}$ along the relation $<_0$. The projection function $\pi: \ig_\iota \to T$ is defined the natural way: each point  in $\ig_\iota$ can be seen as belonging to precisely one $\ig_{\alpha_i}$, so $\pi$ maps it to $P_i$ as $\pi_i$ does.

Let $\chi_0$ be $\bot$ and $\chi_{i}$ be $\lnot \langle 0\rangle \phi_{\alpha_0+\ldots+\alpha_i}$. Clearly, the formula $\lnot \chi_i\land \chi_{i+1}$ defines precisely the copy of $\ig_{\alpha_i}$ within $\ig$. Given $t\in P_i$ we define $\pi^{-1}(t)$ by the formula
$$\psi\land \lnot \chi_i\land \chi_{i+1},$$
where $\psi$ was the formula defining $\pi_i^{-1}(t)$ within $\ig_{\alpha_i}$. Notice that since $\psi$ was constructed as a $[0]$- and $\langle 0\rangle$-free formula its validity in $\ig_{\alpha_i}$ is equivalent to that in its shifted copy within $\ig_\iota$.

Now to finish the proof we need to validate that our construction satisfies property (\ref{LemmaDefinability3}) (property (\ref{LemmaDefinability4}) is implied by (\ref{LemmaDefinability1}) and (\ref{LemmaDefinability3})). However, since our construction essentially is the same as Beklemishev's construction, we refer the reader to \cite[Lemma 9.2]{Be10} for the verification of this property.
\endproof

 %Rather than verifying this, the reader might find it easier to directly verify the following, which is the only property of it that we will require:
% As in \cite[Lemma 9.2]{Be10}, we have that if $\phi\wedge M^+(\phi)$ holds in the root $r$ of some hereditarily linear $J_n$-frame, then $\ig_\iota,\delta_\iota\Vdash\phi$, where $\iota$ is the ordinal obtained in the proof of the lemma.

\section{Proofs of Theorems \ref{TheoremMain} and \ref{TheoremMainGLP}} \label{SectProof}
In this section, we put all the previous results together in order to prove Theorems \ref{TheoremMain} and \ref{TheoremMainGLP}. We begin with a lemma.
\begin{lemma}\label{LemmaZeroConsistent}
Let $\phi$ be a closed $\mathcal{L}_{[<n]}$-formula consistent with $\glp$. Then, there is $k\in\mathbb{N}$ such that
\[\glp \vdash \langle n-1 \rangle^k \top\to\langle 0\rangle \phi.\]
\end{lemma}
\proof
By Ignatiev's theorem, $\phi$ holds in a point along the main axis of $\ig$.
Moreover, such a point can be found in $\ig_\iota$, where
\[\iota < \underbrace{\omega^{\omega^{\iddots^\omega}}}_{n \text{ times}} = \lim_{k\to\infty} \underbrace{\omega^{\omega^{\iddots^{\omega^k}}}}_{n-1 \text{ times}}.\]
This follows from Ignatiev's proof (cf. also the construction in \S\ref{SectDefinability} above).
Note that $\ig_\iota\models [n]\bot$ and thus $\ig_\iota\models [n-1]^k\bot$ for some $k\in\mathbb{N}$; thus,
it follows that $\ig\models \langle n-1\rangle^k \top \to \langle 0 \rangle \phi$, and thus $\glp \vdash \langle n-1 \rangle^k \top \to \langle 0 \rangle \phi$ by completeness.
\endproof

We now have the following theorem. Below, recall that a modal logic is \textit{non-degenerate} if it does not contain $\Box^k \bot$ for any modality.
\begin{theorem}\label{TheoremMainProof}
Let $n\in\mathbb{N}$. The following are equivalent for each $\phi\in\mathcal{L}_{[<n]}$:
\begin{enumerate}
\item\label{TheoremMainProof1} $\phi$ is a theorem of $\glp.3$;
\item\label{TheoremMainProof2} every closed $\mathcal{L}_{[<n]}$-substitution of $\phi$ is a theorem of $\glp_n$;
\item\label{TheoremMainProof3} $\phi$ is a theorem of some non-degenerate, normal extension of $\glp_n$.
\end{enumerate}
%Moreover, if $V = L$, one can also add:
%\begin{enumerate}[resume]
%\item\label{TheoremMainProof4} $\phi$ is valid under the correct-model interpretation.
%\end{enumerate}
\end{theorem}
\proof
%Let us first suppose that $V=L$ holds and prove all equivalences under this assumption. 
Clearly \eqref{TheoremMainProof1} implies \eqref{TheoremMainProof3}. 
%by Theorem \ref{TheoremValidity}. Since $\zfc$ proves that there exist $\Sigma_n$-correct transitive sets, the logic of correct models restricted to the language $\mathcal{L}_{[<n]}$ is normal and non-degenerate (for our purposes, it suffices that it does not contain $[n-1]^k\bot$ for any $k$). Thus, \eqref{TheoremMainProof4} implies \eqref{TheoremMainProof3}.
To see that \eqref{TheoremMainProof3} implies \eqref{TheoremMainProof2}, let $\mathcal{L}$ be a  normal extension of $\glp$ which does not prove $[n-1]^k \bot$ for any $k\in\mathbb{N}$. 
Suppose towards a contradiction that $\mathcal{L}\vdash\phi$ but $\lnot\phi^*$ is consistent with $\glp$, for some closed $\mathcal{L}_{[<n]}$-substitution $\phi^*$ of $\phi$.
Since $\mathcal{L}\vdash\phi$, we have $\mathcal{L}\vdash\phi^*$ and hence $\mathcal{L}\vdash[0]\phi^*$ by normality. Thus we have $\mathcal{L}\vdash\lnot\langle 0 \rangle \lnot \phi^*$, so $\langle 0 \rangle \lnot\phi^*$ is inconsistent with $\mathcal{L}$. However, $\lnot\phi^*$ is a closed $\mathcal{L}_{[<n]}$-formula consistent with $\glp$, so by Lemma \ref{LemmaZeroConsistent} we have
\[\exists k\in\mathbb{N}\Big(\glp \vdash \langle n-1\rangle^k \top \to \langle 0\rangle \lnot \phi^*\Big),\]
contradicting the fact that $\mathcal{L}$ is consistent with $\langle n-1 \rangle^k \top$ for all $k\in\mathbb{N}$.

Finally, to see that \eqref{TheoremMainProof2} implies \eqref{TheoremMainProof1}, we suppose that $\phi$ is consistent with $\glp.3$. Thus, $\phi\wedge M^+(\phi)$ is also consistent with $\glp.3$, since $M^+(\phi)$ is trivially provable in $\glp$. In particular, $\phi\wedge M^+(\phi)$ is consistent with $\ja_n.\lin$ and thus it holds in the root $r$ of some finite hereditarily linear $J$-frame $T$ by Theorem \ref{TheoremCompletenessJlin}. As in \S \ref{SectDefinability} and using Lemma \ref{LemmaDefinability}, we may find ordinal $\iota<\varepsilon_0$ and a closed substitution $\phi^*$ of $\phi$  such that $\ig,\delta_\iota\Vdash \phi^*$, so $\phi^*$ is consistent with $\glp$.
%We have proved the equivalence between \eqref{TheoremMainProof1}--\eqref{TheoremMainProof4} under the assumption that $V = L$. However, observe that each of \eqref{TheoremMainProof1}--\eqref{TheoremMainProof3} is equivalent to a $\Sigma^1_1$ sentence; since $\zfc + V = L$ is $\Sigma^1_1$-conservative over $\zfc$ e.g., by Mostowski's absoluteness theorem, the equivalence between \eqref{TheoremMainProof1}--\eqref{TheoremMainProof3} holds without the hypothesis that $V = L$.
\endproof

Let us draw some conclusions from Theorem \ref{TheoremMainProof}. First, the global version for full $\glp$, i.e., Theorem \ref{TheoremMainGLP} from the introduction, follows immediately.
We also have the following consequence:
\begin{theorem}
$\glp.3$ is the logic of closed substitutions of $\glp$, as well as the largest non-degenerate, normal extension of $\glp$.
\end{theorem}

From Theorem \ref{TheoremMainProof} we can also derive Theorem \ref{TheoremMain} as a schema:

\begin{theorem}\label{TheoremMainProofVL}
Suppose $V = L$ and let $n\in\mathbb{N}$. The following are equivalent for each $\phi\in\mathcal{L}_{[<n]}$:
\begin{enumerate}
\item\label{TheoremMainProofVL1} $\phi$ is a theorem of $\glp.3$; and
\item\label{TheoremMainProof4} $\phi$ is valid under the correct-model interpretation.
\end{enumerate}
\end{theorem}
\proof
First, \eqref{TheoremMainProof1} implies \eqref{TheoremMainProof4} 
by Theorem \ref{TheoremValidity}. Conversely, $\zfc$ proves that there exist $\Sigma_n$-correct transitive sets, the logic of correct models restricted to the language $\mathcal{L}_{[<n]}$ is normal and non-degenerate. Thus, \eqref{TheoremMainProof4} implies \eqref{TheoremMainProofVL1} by Theorem \ref{TheoremMainProof}.
\endproof

\subsection{Discussion}
There are several aspects of the correct-model interpretation which have been chosen rather arbitrarily. For instance, one could modify the interpretation by mapping $[n]\phi$ to a formalization of ``$\phi$ holds in every $\Sigma_n$-correct model satisfying $\kp + V = L$, as well as $\Sigma_1$-Separation.'' For $n \neq 0$, the sentence $V = L$ follows from $\Sigma_n$-correctness; for $n = 0$, it needs to be added separately in order for the proof of Theorem \ref{TheoremValidity} to go through.

One could also modify the interpretation by restricting to submodels of some $V_\kappa$ or $H(\kappa)$. For instance, one could interpret $[n]\phi$ as ``$\phi$ holds in every $M\prec_{\Sigma_{n+1}} H(\omega_1)$.'' Since $\zfc$ can uniformly prove the existence of $\Sigma_n$-correct submodels of $H(\omega_1)$, the corresponding version of Theorem \ref{TheoremMain}  becomes provable when formalized as a single sentence, by the proofs of Theorem \ref{TheoremMainProof} and Theorem \ref{TheoremMainProofVL}.

Theorem \ref{TheoremMain} is no longer true when the hypothesis that $V = L$ is removed. It follows from the results of \cite{AgPac} that the interpetation of $[n]$ described in the previous paragraph has $\mathsf{GL}$ as its unimodal logic if and only if $n$ is even, assuming that sufficiently large cardinals exist. In particular, axiom $.3$ fails under these circumstances. In the following section, we prove that large cardinals are not necessary, and indeed the conclusion of Theorem \ref{TheoremMain} fails as soon as one considers models of set theory which extend $L$ even in a very mild manner.

\section{Proof of Theorem \ref{TheoremVL}}\label{SectTheoremC}
The main theorem of this section is:
\begin{theorem}\label{TheoremCohenRandom}
Suppose that there exist real numbers $c,r\in\mathbb{R}$ such that $c$ is Cohen over $L$ and $r$ is random over $L$. Then, the linearity axiom $.3$ is not sound for the correct-model interpretation, even when restricting to $\mathcal{L}_{[<1]}$-formulas.
\end{theorem}
In particular, the hypotheses of Theorem \ref{TheoremCohenRandom} hold if $\omega_1^L$ is countable, or if Martin's Axiom holds, or if sufficiently large cardinals exist.
We refer the reader to Jech \cite{Je03} for more on generic reals. Here, recall that a real is \textit{Cohen over $L$} if it is of the form $\bigcap g$, where $g\subset \mathbb{C}$ is a filter on the Cohen algebra which has nonempty intersection with every dense subset of $\mathbb{C}$ which belongs to $L$. Similarly, a real is \textit{random over $L$} if it is of the form   $\bigcap g$, where $g\subset \mathbb{M}^L$ is a filter on Solovay's random algebra (as defined in $L$) which has nonempty intersection with every dense subset of $\mathbb{M}^L$ which belongs to $L$.  By a result of Solovay \cite{So70}, a real is Cohen over $L$ if and only if it does not belong to any meager Borel subset of $\mathbb{R}$ with a Borel code in $L$, and a real is random over $L$ if and only if it does not belong to any measure-zero Borel subset of $\mathbb{R}$ with a Borel code in $L$. 

Below, we make use of various standard facts about admissible sets, stability, and $\Sigma_1$-reflection, which we recall for the reader's convenience. 
Further background as well as more detailed proofs of these claims (modulo their straightforward relativizations to real parameters) can be found in Barwise \cite[\S V.7 and \S V.8]{Ba75}.

First, recall that an ordinal $\alpha$ is called \textit{$x$-stable} if $L_\alpha[x] \prec_{\Sigma_1} L[x]$. For $\alpha<\omega_1^{L[x]}$, this can be equivalently stated as
\[L_\alpha[x] \prec_{\Sigma_1} L_{\omega_1^{L[x]}}[x],\]
by the L\"owenheim-Skolem theorem.
Stable ordinals are always admissible, limits of admissible, etc.
We remind the reader of the Shoenfield absoluteness theorem, which asserts that if $\varphi$ is $\Sigma^1_2$ formula without parameters and with only one free variable, then for all $x \in\mathbb{R}$ such that $\varphi(x)$ holds, we have $L[x] \models\varphi(x)$. Moreover, the same holds for $\Sigma_1$ formulas in the language of set theory.

Let $\sigma$ be the least $x$-stable ordinal. Recall that every element of $L_\sigma[x]$ is $\Sigma_1$-definable in $L_\sigma[x]$ with $x$ as the only parameter. The Shoenfield absoluteness theorem thus implies that $L_\sigma[x]$ is $\Sigma_1$-correct. For any $y \in L_\sigma[x]\cap \mathbb{R}$, $\sigma$ will also be $y$-stable. %This is because $y$ is $\Sigma_1$-definable in $L_\sigma[x]$ and thus in $L[x]$, and the function $\gamma\mapsto L_\gamma[y]$ is a total $\Sigma_1$-function, so for any $\Sigma_1$-formula $\psi$, the formula 
%\[\exists \gamma\, L_\gamma[y] \models \psi(y)\]
%if, true, is also true in $L_\sigma[x]$ by stability. Since 
%\[L[y]^{L_\sigma[x]} = L_\sigma[y],\] the claim follows.

Further, if $\sigma$ is an arbitrary $x$-stable ordinal below $\omega_1^{L[x]}$, then one of the following holds: (i) $\sigma$ is a \textit{successor} $x$-stable, in which case $\sigma$ is the least $x'$-stable ordinal relative to some $x' \in L_\sigma[x]$ (namely, one coding $L_{\sigma'}[x]$, where $\sigma'$ is the previous $x$-stable ordinal), or (ii) $\sigma$ is a limit of $x$-stable ordinals. In both cases, it follows from the comment above that $L_\sigma[x]$ is $\Sigma_1$-correct. By similar considerations, if $\sigma$ is an arbitrary  $x$-stable ordinal below $\omega_1^{L[y]}$ and $y \in L_\sigma[x]$, then $\sigma$ is also $y$-stable and thus $L_\sigma[y]$ is also $\Sigma_1$-correct. 

With this background being recalled, we may now proceed to the proof of the theorem.
 %In particular, we remind the reader that to each $\Sigma^1_2$ formula $\phi$ we can uniformly assign a $\Sigma_1$ formula $\phi^*$ such that $(\mathbb{N}, x)\models\phi$ if and only if $L_\sigma[x] \models\phi^*$ whenever $x\in\mathbb{R}$ and  $\sigma$ is an $x$-stable ordinal; and conversely to each $\Sigma_1$ formula $\phi^*$ we can uniformly assign a $\Sigma^1_2$ formula $\phi$ such that for all $x\in\mathbb{R}$ and all $x$-stable ordinals $\sigma$, we have $(\mathbb{N}, x)\models\phi$ if and only if $L_\sigma[x] \models\phi^*$ (this is a version of the well-known Shoenfield absoluteness theorem).
\proof[Proof of Theorem \ref{TheoremCohenRandom}]
We let $\phi_c$ be the sentence ``there are no reals Cohen over $L$'' and $\phi_r$ be ``there are no reals random over $L$.'' We claim that 
\begin{equation}\label{eq3ForcingForm}
[0]\big([0]\phi_c \to \phi_r\big) \vee [0]\big(\phi_r \wedge [0]\phi_r \to \phi_c\big)
\end{equation}
is not valid, i.e., that it fails in some $\Sigma_1$-correct model. Since both Cohen and random forcings have the countable chain condition, $$\omega_1^L=\omega_1^{L[c]}=\omega_1^{L[r]}\le \omega_1^{L[c\oplus r]},$$ where $c\oplus r$ is the Turing join of $c$ and $r$. We construct $\alpha\le \omega_1^L$  that is both $c$ and $r$ stable as the limit of the chain $(\alpha_i)_{i<\omega}$,  where $\alpha_0=0$, $\alpha_{2n+1}$ is a $c$-stable above $\alpha_{2n}$, and $\alpha_{2n+2}$ is an $r$-stable above $\alpha_{2n+1}$. 
Since this construction can be carried out within $L[c\oplus r]$, we may assume $\alpha<\omega_1^{L[c\oplus r]}$.  

Writing $\gamma := \omega_1^{L[c\oplus r]}$,
we have $\Sigma_1$-correct models  $L_\gamma[c\oplus r]$, $L_\alpha[c]$, $L_\alpha[r]$, and $L_\alpha$, and moreover
\begin{align}\label{eqCorrectAlpha}
L_\gamma[c\oplus r]\models  \text{``$L_\alpha[c]$ and $L_\alpha[r]$ are $\Sigma_1$-correct.''}
\end{align}
It suffices to show that \eqref{eq3ForcingForm} fails in $L_\gamma[c\oplus r]$. Suppose otherwise. There are two cases. Suppose first that $[0]([0]\phi_c \to \phi_r)$ holds in $L_\gamma[c\oplus r]$. By \eqref{eqCorrectAlpha}, we have
\[L_\alpha[r] \models [0]\phi_c\to\phi_r.\]
By hypothesis, $r$ does not belong to any measure-zero set with a Borel code in $L$, and thus in $L_\alpha$, so $L_\alpha[r]\models ``r$ is random over $L$.'' It follows that $L_\alpha[r]\models \lnot\phi_r$. Let $M \in L_\alpha[r]$ be transitive and such that $M\prec_{\Sigma_1} L_\alpha[r]$. Towards a contradiction, suppose that $M\models \lnot\phi_c$, and let $x\in M$ be such that $M\models ``x$ is Cohen over $L$.'' It is well known that random forcing adds no Cohen reals (and vice-versa; see Jech \cite{Je03}). Thus, $L[r]\models ``x$ is not Cohen over $L$,'' so $L[r]$ has a code $b$ for a meager Borel set $B$ such that $r \in B^{L[r]}$. 

\begin{claim}\label{ClaimCohen}
The formula ``$x$ is not Cohen over $L$'' is $\Sigma^1_2$.
\end{claim}
\proof
We use Solovay's characterization in terms of Borel codes. Thus, that $x$ is not Cohen over $L$ is equivalent to the existence of $\beta$, $p$, and formulas $\phi$, $\{\phi_i:i\in\mathbb{N}\}$ such that the following hold:
\begin{enumerate}
\item $(\beta, p, \phi, \{ \phi_i : i\in\mathbb{N}\}) \in L$;
\item $\beta$ is an ordinal;
\item $p \in\mathbb{R}$ is a parameter;
\item $\phi$ is a formula in $\Sigma^0_\beta(p)$, as are $\phi_i$ for each $i\in\mathbb{N}$ (i.e., these formulas define Borel sets of rank $\beta$, effectively in $p$);
\item for all $y\in\mathbb{R}$, $\phi(y,p)$ holds if and only if $\phi_i(y,p)$ holds for some $i\in\mathbb{N}$ (i.e., $\phi$ defines the union of the sets defined by the family $\{\phi_i:i\in\mathbb{N}\}$);
\item $\phi(x,p)$ holds (i.e., $x$ belongs to the set coded by $\phi$);
\item for all $y \in\mathbb{R}$ and all $i\in\mathbb{N}$, there exists a rational $\varepsilon>0$ such that for all $z\in\mathbb{R}$, if $\phi_i(z,p)$, then 
$\varepsilon \leq |y - z|.$
\end{enumerate} 
The last condition expresses that the closure of the set defined by $\phi_i$ has empty interior, so that $\phi$ is meager.
\endproof
By the claim, the fact that $x$ is not Cohen over $L$ is $\Sigma_1$ over $L[r]$, so it is true in $M$ by the fact that 
\[M \prec_{\Sigma_1} L_{\alpha}[r]  \prec_{\Sigma_1} L[r],\] which is a contradiction. Thus, we have $M\models \phi_c$. Since $M$ was arbitrary, we have $L_\alpha[r]\models [0]\phi_c$, which is a contradiction.

Suppose now that $[0](\phi_r\wedge [0]\phi_r \to \phi_c)$ holds in $L_\gamma[c\oplus r]$. Using the formula \eqref{eqCorrectAlpha} one more, we see that
\begin{equation}\label{eqGLP3Cohen}
L_\alpha[c]\models \phi_r \wedge [0]\phi_r \to \phi_c.
\end{equation}
As mentioned earlier, Cohen forcing adds no random reals. Moreover, the partial order for adding Cohen reals is simply the set of finite partial functions from $\mathbb{N}$ to $\{0,1\}$ ordered by extension, so it is recursive and thus belongs to $L_{\omega+1}$ and in particular to $L_\alpha$. It follows that $c$ is generic over $L_\alpha$, so $L_\alpha[c]$ has no reals which are random over $L$. The same argument from before, using Claim \ref{ClaimCohen}, shows that $L_\alpha[c]\models$``there are no reals random over $L$,'' i.e., $L_\alpha[c]\models\phi_r$. Thus, \eqref{eqGLP3Cohen} can be strengthened to 
\[L_\alpha[c]\models [0]\phi_r \to \phi_c.\]
As mentioned before, we have $L_\alpha[c]\models ``c$ is Cohen over $L$,''  and thus $L_\alpha[c]\models\lnot [0]\phi_r$.
\begin{claim}\label{ClaimRandom}
The formula ``$x$ is not random over $L$'' is $\Sigma^1_2$.
\end{claim}
\proof
Similar to Claim \ref{ClaimCohen}.
\endproof
However, the same argument from before shows, using Claim \ref{ClaimRandom},  that $L_\alpha[c]\models [0]\phi_r$, which is a contradiction once more.
\endproof

\subsection*{Acknowledgements}
The first author was partially supported by FWF grants ESP-3N and P-36837.  The work of the second author has been funded by the FWO grant
G0F8421N.

\bibliographystyle{abbrv}
\bibliography{References}

\end{document}